\documentclass[12pt]{article}
\usepackage{amsmath,amssymb,amsthm,amscd,graphics,epsfig}
\usepackage[english]{babel}
\usepackage{setspace}

\newtheorem{theorem}{Theorem}
\newtheorem{lemma}[theorem]{Lemma}
\newtheorem{proposition}[theorem]{Proposition}

\newtheorem{definition}[theorem]{Definition}
\newtheorem{remark}[theorem]{Remark}

\def\CC{\mathbb{C}}
\def\DD{\mathbb{D}}
\def\HH{\mathbb{H}}

\def\RR{\mathbb{R}}
\def\SS{\mathbb{S}}

\def\aaa{\mathcal{A}}
\def\cc{\mathcal{C}}
\def\dd{\mathcal{D}}
\def\ii{\mathcal{I}}
\def\ll{\mathcal{L}}
\def\oo{\mathcal{O}}

\def\ZZ{\mathbb{Z}}

\def\Un{\mathop{\text{\bf 1}}\nolimits}

\renewenvironment{proof}[1][.]{%
\bigskip\noindent{\bf Proof#1 }}{%
\hfill$\blacksquare$\bigskip}

\begin{document}

\title{Eigenfunctions  of the Laplacian and  associated
Ruelle operator }
\author{ A. O. Lopes\footnote{ Instituto de Matem\'atica, UFRGS, Porto
Alegre 91501-970, Brasil. Partially supported by CNPq, Fapergs,
Instituto do Milenio, benefici\'ario de aux\'{\i}lio financeiro
CAPES - Brasil.}\and Ph. Thieullen\footnote{ Institut de
Math\'ematiques, CNRS, Universit\'e  Bordeaux 1, 341 cours de
la lib\'eration, F-33405 Talence, France.} } \maketitle

\begin{abstract}
Let $\Gamma$ be a co-compact Fuchsian group of isometries on the
Poincar\'e disk $\DD$ and  $\Delta$ the corresponding hyperbolic
Laplace operator. Any  smooth eigenfunction $f$ of $\Delta$,
equivariant by $\Gamma$ with real eigenvalue $\lambda=-s(1-s)$,
where $s=\frac{1}{2}+ it$, admits an integral representation by a
distribution $\dd_{f,s}$ (the Helgason distribution) which is
equivariant by $\Gamma$ and supported at infinity
$\partial\DD=\SS^1$. The geodesic flow on the compact surface
$\DD/\Gamma$ is conjugate to a suspension over a natural extension
of a piecewise analytic map $T:\SS^1\to\SS^1$, the so-called
Bowen-Series transformation. Let $\ll_s$ be the complex Ruelle
transfer operator associated to the jacobian $-s\ln |T'|$. M.
Pollicott showed that $\dd_{f,s}$ is an eigenfunction of the dual
operator $\ll_s^*$ for the eigenvalue 1. Here we show the existence of
a (nonzero) piecewise real analytic  eigenfunction $\psi_{f,s}$ of
$\ll_s$ for the eigenvalue 1, given by an integral formula
\[
\psi_{f,s} (\xi)=\int \frac{J(\xi,\eta)}{|\xi-\eta|^{2s}} \,
\dd_{f,s} (d\eta),
\]
\noindent where $J(\xi,\eta)$ is a $\{0,1\}$-valued piecewise
constant function whose definition depends upon the geometry of
the Dirichlet fundamental domain representing the surface
$\DD/\Gamma$.
\end{abstract}

\section{Introduction}

Consider the Laplace operator $\Delta$ defined by
\[
\Delta=y^2\left(\,\frac{\partial^2 }{\partial x^2}+
\frac{\partial^2 }{\partial y^2}\,\right),\] on the Lobatchevskii
upper half-plane $\HH=\{w=x+iy\in\CC; y>0\}$, equipped with the
hyperbolic metric $ds_{\HH}=\frac{|dw|}{y}$, and the eigenvalue
problem
\[
\Delta f=-s(1-s)f,
\]
\noindent where $s$  is of the form $s=\frac{1}{2}+it$, with $t$
is real. We shall also consider the same corresponding Laplace
operator
\[
\Delta=\frac{1}{4}(1-|z|^2)^2\left(\frac{\partial^2}{\partial
x^2}+\frac{\partial^2}{\partial y^2}\right),\] and eigenvalue
problem
\[
\Delta f=-s(1-s)f,
\]
defined on the Poincar\'e disk $\DD=\{z=x+yi\in\CC; |z|<1\}$,
equipped with the metric $ds_{\DD}=2\frac{|dz|}{1-|z|^2}$.

Helgason showed in  \cite{He1} and \cite{He} that any
eigenfunction $f$ associated to this eigenvalue problem can be
obtained by means of a generalized Poisson representation
\[
\left\{\begin{array}{ll}
\displaystyle
f(w)=\int_{-\infty}^\infty \left(\frac{(1+t^2)y}{(x-t)^2 +y^2}\right)^{s}\,\dd_{f,s}^{\HH}(t),&\text{for } w\in\HH,\\
\\
\text{or} \\
\displaystyle f(z)=\int_{\partial\,
\DD}\left(\frac{1-|z|^2}{|z-\xi|^{2}}\right)^s\,\dd_{f,s}^{\DD}(\xi),&\text{for
} z\in\DD,
\end{array}\right.
\]
\noindent where $\dd_{f,s}^{\DD}$ or $\dd_{f,s}^{\HH}$ are
analytic distributions called from now on {\it Helgason's
distributions}. We have used the canonical isometry between
$z\in\DD$ and $w\in\HH$, namely $w=i\frac{1-z}{1+z}$ or
$z=\frac{i-w}{i+w}$. The hyperbolic metric is given in $\HH$ and
in $\DD$ by
\[
ds^2_\HH=\frac{dx^2+dy^2}{y^2},\quad
ds^2_\DD=\frac{4(dx^2+dy^2)}{(1-|z|^2)^2}.
\]

We shall be interested in a more restricted problem, where the
eigenfunction $f$ is also automorphic with respect to a co-compact
Fuchsian group $\Gamma$,  i. e., a discrete subgroup of the group
of M\"obius transformations (see \cite{Pa}, \cite{Su},
\cite{BeMa}) with compact fundamental domain. It is known that the
eigenvalues $\lambda=s(1-s)=\frac{1}{4}+t^2$ form a discrete set
of positive real numbers with finite multiplicity and accumulating
at $+\infty$ (see \cite{I}).

M. Pollicott showed  \cite{Po} that the Helgason's distribution can
be seen as a generalized eigenmeasure of the dual  complex Ruelle
transfer operator associated to a subshift of finite type defined
at infinity. Let $T_L$ be the left Bowen-Series transformation
that acts on the boundary $\SS^1=\partial \DD$ and is  associated
to a particular set of generators of $\Gamma$. The precise
definition of $T_L$ has been given  in \cite{BoSe1}, \cite{Ser2},
\cite{Ser3}, \cite{Ser1}, and  more geometrical descriptions have
then been given in \cite{AdFl} and \cite{Mo}. Specific examples of
the Bowen-Series transformation have been studied in \cite{Ma} and
\cite {BeKeSe} for the modular surface and in \cite{BLT} for a
symmetric compact fundamental domain of genus two. The map $T_L$
is known to be piecewise $\Gamma$-M\"obius constant, Markovian
with respect to a partition $\{I_k^L\}$ of intervals of $\SS^1$,
on which the restriction of $T_L$ is constant and equal to an
element $\gamma_k$ of $\Gamma$, transitive and orbit equivalent to
$\Gamma$. Let $\ll_s^L$ be the {\it complex Ruelle transfer
operator} associated to the map $T_L$ and the potential
$A_L=-s\ln|T'_L|$, namely
\[
(\ll_s^L\psi)(\xi') =\sum_{T_L(\xi)=\xi'}e^{A_L(\xi)}\psi(\xi)
=\sum_{T_L(\xi)=\xi'}\frac{\psi(\xi)}{|T'_L(\xi)|^s},
\]
\noindent where the summation is taken over all  pre-images $\xi$
of $\xi'$ under $T_L$. Here  $T'_L$ denotes the Jacobian of $T_L$
with respect to the canonical Lebesgue measure on $\SS^1$. In the
case of an automorphic eigenfunction $f$ of $\Delta$, Pollicott
showed that the corresponding Helgason distribution $\dd_{f,s}$
satisfies the dual functional equation
\[
(\ll_s^L)^*(\dd_{f,s})=\dd_{f,s}
\]
\noindent or, according to Pollicott's terminology, the parameter
$s$ is a {\it (dual) Perron-Frobenius value}, that is, 1 is an
eigenvalue for the dual Ruelle transfer operator.

Although suggested in \cite{Po}, it is not clear whether $s$ could
be a  {\it Perron-Frobenius value}, that is, whether 1 could also
be an eigenvalue for $\ll_s^L$, not only for $(\ll_s^L)^*$. Our
goal in this paper is to show that this is actually the case.

The three main ingredients we use are the following;
\begin{itemize}
\item Otal's proof of Helgason's distribution in \cite{Ot}, giving
more precise information on $\dd_{f,s}$ and enabling us to
integrate piecewise $\cc^1$ test functions, instead of real
analytic globally defined test functions; \item a more careful
reading of \cite{AdFl}, \cite{Mo}, \cite{BoSe1}, and  \cite{Ser1},
or a careful study of a particular example in \cite{LT}, which
enables us to construct a piecewise $\Gamma$-M\"obius baker
transformation (``arithmetically'' conjugate to the geodesic
billiard); \item the existence of a kernel that we introduced in
\cite{BLT}, which enables us to permute past and future
coordinates and transfer a dual eigendistribution to a piecewise
real analytic eigenfunction. Haydn (in \cite{Hay}) has introduced
a similar kernel in a more abstract setting, without geometric
considerations.
\end{itemize}

More precisely, we prove the following theorem

\begin{theorem} \label{label0}
Let $\Gamma$ be a co-compact Fuchsian group of the hyperbolic disk
$\DD$ and $\Delta$  the corresponding hyperbolic Laplace operator.
Let $\lambda=s(1-s)$, with $s=\frac{1}{2}+it$, and let $f$ be an
eigenfunction of $-\Delta$, automorphic with respect to $\Gamma$,
that is, $\Delta f=-\lambda f$ and $f\circ\gamma=f$, for every
$\gamma\in\Gamma$. Then there exists a (nonzero) piecewise real
analytic eigenfunction $\psi_{f,s}$  on $\SS^1$ that is a solution
of the functional equation
\[
\ll_s^L(\psi_{f,s})=\psi_{f,s},
\]
\noindent where $\ll^L_s$ is the complex Ruelle transfer operator associated to the left Bowen-Series transformation $T_L:\SS^1\to\SS^1$ and the potential $A_L=-s\ln|T'_L|$.

Moreover, $\psi_{f,s}$ admits an integral representation via
Helgason's distribution $\dd_{f,s}^{\DD}$, representing $f$ at
infinity, and a geometric positive kernel $k(\xi,\eta)$ defined on
a finite set of disjoint  rectangles $\cup_k I^L_k\times
Q_k^R\subset \SS^1\times\SS^1$, namely,
\[
\psi_{f,s}(\xi)=\int_{Q^R_k}
k^s(\xi,\eta)\,\dd_{f,s}^{\DD}(\eta)=\int_{Q^R_k}\frac{1}{|\xi-\eta|^{2s}}\,\dd_{f,s}^{\DD}(\eta),
\qquad
\]
\noindent for every $ \xi\in I^L_k$, where $I^L_k$ and $Q^R_k$ are
intervals of $\SS^1$ on which disjoint closure, and $\{I^L_k\}_k$
is a partition of $\SS^1$ where $T_L$ is injective, Markovian  and
piecewise $\Gamma$-M\"obius constant.
\end{theorem}

Lewis \cite{Le} and, later, Lewis and Zagier \cite{LeZa}, started
a different approach to understand Maass wave forms. They where
able to identify in a bijective way Maass wave forms of
$PSL(2,\ZZ)$ and solutions of a functional equation with 3 terms
closely related to Mayer's transfer operator. Their setting is
strongly dependent of the modular group.  Our theorem \ref{label0}
may be viewed as part of their program for co-compact Fuchsian
groups. The Helgason distribution has been used by S. Zelditch in
\cite{Ze} to generalise microlocal analysis on hyperbolic
surfaces, by L. Flaminio and G. Forni in \cite{FF}, to study
invariant distributions by the horocycle flow, and by N.
Anantharaman and S. Zelditch in \cite{AZ}, to understand the
``Quantum Unique Ergodicity Conjecture''.

\section{Preliminary results}
\vspace{0.3cm}

Let $\Gamma$ be a co-compact Fuchsian group of the Poincar\'e disk
$\DD$. We denote by $d(w,z)$ the hyperbolic distance between two
points of $\DD$, given by the Riemannian metric
$ds^2=4(dx^2+dy^2)/(1-|z|^2)^2$. Let $M=\DD/\Gamma$ be the
associated compact Riemann surface,  $N=T^1M$ the unit tangent
bundle, and  $\Delta$ the Laplace operator on $M$. Let
$f:M\to\RR$ be an eigenfunction of $-\Delta$ or, in other words, a
$\Gamma$-automorphic function $f:\DD\to\RR$ satisfying $\Delta
f=-s(1-s)f$ for the eigenvalue $\lambda=s(1-s)>\frac{1}{4}$ and
such that  $f\circ\gamma=f$, for every $\gamma\in\Gamma$. We know
that $f$ is $\cc^\infty$ and uniformly bounded on $\DD$. Thanks to
Helgason's representation theorem, $f$ can be represented as a
superposition of horocycle waves, given by the {\it Poisson
kernel}
\[
P(z,\xi):=e^{b_\xi(\oo,z)}=\frac{1-|z|^2}{|z-\xi|^2},
\]
\noindent where $b_\xi(w,z)$ is the the {\it Busemann cocycle}
between two points $w$ and $z$ inside the Poincar\'e disk,
observed from a point at infinity $\xi\in\SS^1$, defined by
\[
b_\xi(w,z):=\text{ ``$d(w,\xi)-d(z,\xi)$''
}=\lim_{t\rightarrow\,\, \xi}d(w,t)-d(z,t),
\]
\noindent where the limit is uniform in $t\rightarrow \xi$ in any
hyperbolic cone at $\xi$. Helgason's theorem states that
\[
f(z)=\int_\DD P^s(z,\xi)\,\dd_{f,s}(\xi)=\langle \dd_{f,s},P^s(z,.)\rangle
\]
for some analytic distribution $\dd_{f,s}$ acting on real analytic
functions  on $\SS^1$. Unfortunately, Helgason's work is too
general and is valid for any eigenfunction not necessarily
equivariant by a group. For bounded $\cc^2$ functions $f$, Otal
\cite{Ot} has shown that the distribution $\dd_{f,s}$ has stronger
properties and can be defined in a simpler manner.

We first recall some standard notations in hyperbolic geometry. We
call $d(z,z_0)$ the hyperbolic distance between two points: for
instance, the distance from the origin is given by
$d(\oo,\tanh(\frac{r}{2})e^{i\theta})=r$.  Let $\cc(\oo,r)$ denote
the set of points in $\DD$ at  hyperbolic distance $r$ from the
origin,
\[
\cc(\oo,r)=\{z\in\DD; \,|z|=\tanh(\textstyle\frac{r}{2})\}
\]
and, more generally, given for any interval $I$ at infinity and
any point $z_0\in\DD$, let $\cc(z_0,r,I)$ denote the angular arc
at the hyperbolic distance $r$ from $z_0$ delimited at infinity by
$I$, that is,
\[
\cc(z_0,r,I)=\{z\in\DD;\,z\in[[z_0,\xi]]\text{ for some }\xi\in I
\,\text{and} \,d(z,z_0)=r\,\},
\]
where $[[z_0,\xi]]$ denotes the geodesic ray from $z_0$ to the
point $\xi$ at infinity. Let $\frac{\partial}{\partial
n}=\frac{\partial}{\partial r}$ denote the exterior normal
derivative to $\cc(\oo,r)$ and $|dz|_{\DD}=\sinh(r)\,d\theta$the
hyperbolic arc length on $\cc(\oo,r)$.

\begin{theorem} {\rm (\cite{Ot} )} \label{theorem:Otal} Let $f$ be a bounded  $\cc^2$ eigenfunction satisfying $\Delta f=-s(1-s) f$.
Then:
\begin{enumerate}
\item There exists a continuous linear functional $\dd_{f,s}$
acting on $\cc^1$ functions of $\SS^1$,defined by
\[
\int\psi(\xi)\,\dd_{f,s}(\xi):=\lim_{r\rightarrow+\infty}\frac{1}{c(s)}\int_{\cc(\oo,r)}\psi(z)e^{-sr}
\left(\frac{\partial f}{\partial n}+sf\right)\,|dz|_{\DD},
\]
where $c(s)$ is a non zero normalizing constant such that $\langle
\dd_{f,s}, \boldsymbol{1}\rangle=f(0)$, and $\psi(z)$ is any
$\cc^1$ extension of $\psi(\xi)$ to a neighborhood of $\SS^1$.
\item $\dd_{f,s}$ represents $f$ in the following sense:
\[
f(z)=\int \big[P(z,\xi)\big]^s\,\dd_{f,s}(\xi),\quad\forall\ z\in\DD.
\]
\noindent $\dd_{f,s}$ is unique and is called the Helgason
distribution of $f$. \item For all $0\leq \alpha \leq 2\pi$, the
following limit exists:
\[
\tilde\dd_{f,s}(\alpha):=\lim_{r\to+\infty} \frac{1}{c(s)} \int_0^\alpha e^{-sr} \left(\frac{\partial f}{\partial n}+sf\right)\Big(\tanh(\frac{r}{2})e^{i\theta}\Big) \sinh(r)d\theta.
\]
\noindent The convergence is uniform in $\alpha\in[0,2\pi]$ and
$\tilde\dd_{f,s}(0)=0$. \item $\tilde\dd_{f,s}$ can be extended to
$\RR$ as  a $\frac{1}{2}$-H\"older continuous function satisfying:
\begin{enumerate}
\item $\tilde\dd_{f,s}(\theta+2\pi)=\tilde\dd_{f,s}(\theta)+f(0)$,
for every $\theta\in\RR$, \item for any $\cc^1$ function
$\psi:\SS^1\to\CC$, denoting $\tilde \psi(\theta)=\psi(\exp
i\theta)$,
\[
\int \psi(\xi)\,\dd_{f,s}(\xi)=\tilde\psi(0)f(0)-\int_0^{2\pi}\frac{\partial\tilde\psi}{\partial\theta}\tilde\dd_{f,s}(\theta)\,d\theta.
\]
\end{enumerate}
\end{enumerate}
\end{theorem}

\noindent Using similar technical tools as Otal, one can prove the
following extension of $\dd_{f,s}$ on piecewise $\cc^1$ functions,
that is, on functions not necessarily continuous but which admit a
$\cc^1$ extension on each interval $[\xi_k,\xi_{k+1}]$ of some
finite and ordered subdivision $\{\xi_0,\xi_1,\dots,\xi_{r-1}\}$
of $\SS^1$.

\begin{proposition}\label{proposition:extension} Let $f$ and $\dd_{f,s}$ be as in Theorem
\ref{theorem:Otal}.
\begin{enumerate}
\item For any interval $I\subset\SS^1$ and any function
$\psi:I\to\CC$, which is $\cc^1$ on the closure of $I$ and null
outside $I$, the following limit exists:
\[
\int \psi(\xi)\,\dd_{f,s}(\xi) := \frac{1}{c(s)}\lim_{r\rightarrow+\infty}
\int_{\cc(\oo,r,I)}\psi(z)e^{-sr}\left(\frac{\partial f}{\partial n}+sf\right)\,|dz|_{\DD}
\]
\noindent where again $\psi(z)$ is any $\cc^1$ extension of
$\psi(\xi)$ to a neighborhood of $\SS^1$. \item For any $0\leq
\alpha<\beta \leq 2\pi$ and  any $\cc^1$ function  $\psi$ on the
interval $I=[\exp(i\alpha),\exp(i\beta)]$,
\[
\int \psi(\xi)\,\dd_{f,s}(\xi) =
\tilde\psi(\beta)\tilde\dd_{f,s}(\beta) -
\tilde\psi(\alpha)\tilde\dd_{f,s} (\alpha) -
\int_\alpha^\beta\frac{\partial\tilde\psi}{\partial\theta}\tilde\dd_{f,s}(\theta)\,d\theta,
\]
\noindent where $\tilde\dd_{f,s}$ and $\tilde\psi(\theta)$ have
been defined in Theorem \ref{theorem:Otal}.
\end{enumerate}
\end{proposition}

\begin{proof}
Given $\alpha\in[0,2\pi]$, let $I=\{e^{i\theta}\mid 0\leq
\theta\leq \alpha\}$ be an interval in $S^1$, and  $\psi$ a
$\mathcal{C}^1$ function defined on a neighborhood of $\SS^1$.
Denote $\tilde \psi(r,\theta)=\psi(\tanh(\frac{r}{2})e^{i\theta})$
and $K(r,\theta)=e^{-sr}\Big(\frac{\partial f}{\partial
n}+sf\Big)\Big(\tanh(\frac{r}{2}e^{i\theta}\Big)\sinh(r)$. Then
\begin{align*}
\frac{1}{c(s)}&\int_{\cc(\oo,r,I)}\psi(z)e^{-sr}\left(\frac{\partial f}{\partial n}+sf\right)\,|dz|_{\DD} \\
&=\int_0^\alpha \tilde\psi(r,\beta)K(r,\beta) \,d\beta \\
&=\int_0^\alpha\Big[\tilde\psi(r,\alpha)+\int_\beta^\alpha-\frac{\partial\tilde\psi}{\partial \theta}(r,\theta)\,d\theta\Big]K(r,\beta)\,d\beta \\
&= \tilde\psi(r,\alpha)\int_0^\alpha K(r,\beta)\,d\beta - \int_0^\alpha \frac{\partial\tilde\psi}{\partial \theta}(r,\theta) \Big[\int_0^\theta K(r,\beta)\,d\beta\Big]\,d\theta.
\end{align*}
\noindent Since $\int_0^\alpha K(r,\beta)\,d\beta \rightarrow
\tilde\dd_{f,s}(\alpha)$ uniformly in $\alpha\in[0,2\pi]$, the
left-hand side of the previous equality converges to
\[
\int \psi(\xi)\Un_{\{\xi\in I\}}\,\dd_{f,s}(\xi) = \tilde\psi(\alpha)\tilde\dd_{f,s}(\alpha) - \int_0^\alpha\frac{\partial\tilde\psi}{\partial\theta}(\theta)\tilde\dd_{f,s}(\theta)\,d\theta.
\]
\noindent The second part of the proposition follows subtracting
such an expression from another one, such as
\[
\int\psi(\xi)\Un_{\{\xi=e^{i\theta} ;\,0\leq \theta \leq
\beta\}}\tilde\dd_{f,s}(\xi) -
\int\psi(\xi)\Un_{\{\xi=e^{i\theta};\, 0\leq \theta \leq
\alpha\}}\tilde\dd_{f,s}(\xi).
\]
\end{proof}

If, in addition, we assume that $f$ is equivariant with respect to
a co-compact Fuchsian group $\Gamma$, Pollicott observed in
\cite{Po} that $\dd_{f,s}$, acting on real analytic functions, is
equivariant by $\Gamma$, that is, satisfies
$\gamma^*(\dd_{f,s})(\xi)=|\gamma'(\xi)|^s\dd_{f,s}(\xi)$, for all
$\gamma\in\Gamma$. Because Otal's construction is more precise and
implies that Helgason's distribution also acts on piecewise
$\cc^1$ functions, the above equivariance property can be improved
in the following way.

\begin{proposition}\label{label1}
Let $f:\DD\to\RR$ be a $\cc^2$ function, $I\subset \SS^1$  an
interval and $\psi:I\to\CC$  a  $\cc^1$ function on the closure of
$I$. If $f$ satisfies $f\circ\gamma=f$, for some
$\gamma\in\Gamma$, ($f$ is not necessarily automorphic), then
\[
\langle \dd_{f,s}, \frac{\psi\circ\gamma^{-1}}{|\gamma'\circ\gamma^{-1}|^s} \Un_{\gamma(I)}\rangle  = \langle \dd_{f,s},\psi\Un_{I}\rangle.
\]
\end{proposition}

The main difficulty here is to transfer the equivariance property
$f\circ\gamma=f$ to an equivalent property for the extension of
$\dd_{f,s}$ to piecewise $\cc^1$ functions. If $I=\SS^1$ and
$\psi$ is real analytic, then, by uniqueness of the
representation, Proposition \ref{label1} is easily proved. It
seems that just knowing the fact that $\dd_{f,s}$ is the
derivative of some H\"older function is not enough to reach a
conclusion. The following proof uses Otal's approach and,
essentially, the extension of $\dd_{f,s}$ described in Part 1 of
Proposition \ref{proposition:extension}.

\begin{proof}[ of Proposition \ref{label1}.]
First we prove the proposition for $\psi=1$. Let
$g(z)=\exp(-sd(\oo,z))$. By definition of $\dd_{f,s}$, we obtain
\begin{align*}
\int \Un_I(\xi)\,\dd_{f,s}(\xi) &= \lim_{r\to+\infty}\frac{1}{c(s)}\int_{\cc(\oo,r',I)}\left(g\frac{\partial f}{\partial n}-f\frac{\partial g}{\partial n}\right)\,|dz|_{\DD}\\
&= \lim_{r\to+\infty}\frac{1}{c(s)}
\int_{\cc(\oo',r',\gamma(I))}\left(g'\frac{\partial f}{\partial
n}- f\frac{\partial g'}{\partial n}\right)\,|dz|_{\DD},
\end{align*}
where $r'=r+d(\oo,\oo')$, $\oo'=\gamma(\oo)$ and
$g'=g\circ\gamma^{-1}$. Notice that the domain bounded by the
circle $\cc(\oo',r')$ contains the circle $\cc(\oo,r)$. Let
$\overline{PQ}$ be the positively oriented arc
$\cc(\oo,r,\gamma(I))$ and  $\overline{P'Q'}$ be the arc
$\cc(\oo',r',\gamma(I))$. Then the two geodesic segments
$[[P,P']]$ and $[[Q,Q']]$ belong to the annulus $r\leq
d(z,\oo)\leq r+2d(\oo,\oo')$ and their length is uniformly
bounded.

We now use Green's formula to compute the right hand side of the
above expression. Let $\Omega$ denote the domain delimited by
$P,P',Q',Q$ using the corresponding arcs and geodesic segments,
and let  $dv=\sinh(r)\,dr\,d\theta$ be the hyperbolic volume
element. We obtain
\begin{align*}
\int_{\overline{P'Q'}}\left(g'\frac{\partial f}{\partial n}-f\frac{\partial g'}{\partial n}\right)\,
|dz|_{\DD}&=\int_{\overline{PQ}}\left(g'\frac{\partial f}{\partial n}-f\frac{\partial g'}{\partial n}\right)\,|dz|_{\DD}\\
&\quad-\int_{[[P,P']]}\cdots\,|dz|_{\DD}-\int_{[[Q',Q]]}\cdots\,|dz|_{\DD}\\
&\quad+\iint_\Omega\left(g'\Delta f-f\Delta g'\right)\,dv.
\end{align*}
When $r$ tends to infinity, the last three terms at the right-hand
side tend to $0$, since along the geodesic segments $[P,P']$ and
$[Q,Q']$, the gradient $\nabla g'$ is uniformly bounded by
$\exp(-\frac{1}{2}r)$ and
\[
g'\Delta f-f\Delta g'=sg'f\sinh(d(z,\oo'))^{-2}
\quad\textrm{and}\quad
\frac{\partial}{\partial n} g'+sg'
\]
are uniformly bounded  by a constant times  $\exp(-\frac{5}{2}r)$
in the domain $\Omega$, for the first expression, and by a
constant times $\exp(-\frac{3}{2}r)$ on $\cc(\oo,r)$, for the
second expression. It follows that
\begin{align*}
\int \Un_I(\xi)\,\dd_{f,s}(\xi) &= \lim_{r\to+\infty}\frac{1}{c(s)}\int_{\cc(\oo,r,\gamma(I))}g' \left(\frac{\partial f}{\partial n}+sf\right)\,|dz|_{\DD}\\
&=\lim_{r\to+\infty}\frac{1}{c(s)}\int_{\cc(\oo,r,\gamma(I))}\big[\psi(z)\big]^s
e^{-sr} \left(\frac{\partial f}{\partial n}+sf\right)\,|dz|_{\DD},
\end{align*}
where $\psi(z)=\exp \left(d(\oo,z)-d(\oo,\gamma^{-1}(z))\right)$.
Now we observe  that
\[
\left\{\begin{array}{ll}
\psi(z) = \exp s\left(d(\oo,z)-d(\gamma(\oo),z)\right), &\quad\textrm{for $z\in\DD$},\\
\psi(\xi) = \exp b_\xi(\oo,\gamma(\oo))= |\gamma'\circ\gamma^{-1}(\xi)|^{-1},
&\quad\textrm{for $\xi\in\partial\DD$},
\end{array}\right.
\]
\noindent actually coincides with a real analytic function
$\Psi(z)$ defined in a neighborhood of $\SS^1$, given explicitly
by
\[
\Psi(z)=\left(
\frac{(1+|z|)^2}{(1+|\gamma^{-1}(z)|)^2|\gamma'\circ\gamma^{-1}(z)|}
\right)^s.
\]
\noindent Thus we proved  proved that
\[
\int \Un_I(\xi)\,\dd_{f,s}(\xi)=\int \frac{\Un_{\gamma(I)}(\xi)}{|\gamma'\circ\gamma^{-1}(\xi)|^s}\,\dd_{f,s}(\xi).
\]

\medskip
\noindent Now we  prove the general case. We use the same notation
for the lifting $\gamma:\RR\mapsto\RR$ of a M\"obius
transformation $\gamma:\SS^1\mapsto\SS^1$. The lifting satisfies
$\gamma(\alpha+2\pi)=\gamma(\alpha)+2\pi$,
$\exp(i\gamma(\alpha))=\gamma(\exp(i\alpha))$ and
$\gamma'(\alpha)=|\gamma'(\alpha)|$, for all $\alpha\in\RR$. Using
Proposition \ref{proposition:extension}, we obtain
\begin{align*}
\tilde\dd_{f,s}(\beta)&-\tilde\dd_{f,s}(\alpha)\\
&= \frac{\tilde\dd_{f,s}\circ\gamma(\beta)}{\gamma'(\beta)^s} - \frac{\tilde\dd_{f,s}\circ\gamma(\alpha)}{\gamma'(\alpha)^s} - \int_{\gamma(\alpha)}^{\gamma(\beta)}\frac{\partial}{\partial\theta} \left(\frac{1}{(\gamma'\circ\gamma^{-1}(\theta))^s}\right) \tilde\dd_{f,s}(\theta)\,d\theta.
\end{align*}
\noindent For any  $\cc^1$ function $\psi(\xi)$ defined on $I$, we
denote $\tilde\psi(\theta)=\psi(\exp i\theta))$, and obtain

\begin{align*}
&LHS := \int \psi(\xi)\Un_I(\xi) \,\dd_{f,s}(\xi) \\
&= \tilde\psi(\beta)\tilde\dd_{f,s}(\beta) - \tilde\psi(\alpha)\tilde\dd_{f,s}(\alpha) - \int_{\alpha}^{\beta}  \frac{\partial\tilde\psi}{\partial\theta}
\tilde\dd_{f,s}(\theta) \,d\theta \\
&= \tilde\psi(\beta)\tilde\dd_{f,s}(\beta) -  \tilde\psi(\alpha)\tilde\dd_{f,s}(\alpha) - \int_{\gamma(\alpha)}^{\gamma(\beta)} \frac{\partial}{\partial\theta} \left(\tilde\psi\circ\gamma^{-1}(\theta)\right) \tilde\dd_{f,s}(\gamma^{-1}\theta) \,d\theta \\
&= \tilde\psi(\beta)\left( \tilde\dd_{f,s}(\beta) - \tilde\dd_{f,s}(\alpha)\right)
 - \int_{\gamma(\alpha)}^{\gamma(\beta)} \frac{\partial\tilde\psi(\gamma^{-1}\theta)}{\partial\theta} \left(\tilde\dd_{f,s}(\gamma^{-1}\theta) - \tilde\dd_{f,s}(\alpha)\right) \,d\theta.
\end{align*}
\noindent We now use the above equivariance and replacee both
$\tilde\dd_{f,s}(\beta) - \tilde\dd_{f,s}(\alpha) $ and
$\tilde\dd_{f,s}(\gamma^{-1}\theta) - \tilde\dd_{f,s}(\alpha)$ by
the corresponding formula involving
$\tilde\dd_{f,s}\circ\gamma(\beta)$,
$\tilde\dd_{f,s}\circ\gamma(\alpha)$, $\tilde\dd_{f,s}(\theta)$.
Thus
\begin{align*}
&LHS \\
&= \frac{\tilde\psi(\beta)\tilde\dd_{f,s}\circ\gamma(\beta)}{\gamma'(\beta)^s} - \frac{\tilde\psi(\alpha)\tilde\dd_{f,s}\circ\gamma(\alpha)}{\gamma'(\alpha)^s} - \int_{\gamma(\alpha)}^{\gamma(\beta)}\frac{\partial}{\partial\theta} \left(\frac{\tilde\psi(\gamma^{-1}\theta)}{\gamma'(\gamma^{-1}\theta)^s}\right) \tilde\dd_{f,s}(\theta)\,d\theta \\
&=\int \frac{\psi\circ\gamma^{-1}(\xi)}{|\gamma'\circ\gamma^{-1}(\xi)|^s}\Un_{\gamma(I)} \,\dd_{f,s}(\xi).
\end{align*}
\end{proof}

Following \cite{AdFl}, \cite{BoSe1}, \cite{Ser2}, \cite{Ser3},
\cite{Ser1} and  \cite{Mo} for the general case and \cite{LT} for
a specific example we recall the definition of the left $T_L$ and
right $T_R$ Bowen-Series transformation. The hyperbolic surface we
are interested in is given by the quotient of the hyperbolic disk
$\DD$ by a co-compact Fuchsian group $\Gamma$. Given a point
$\oo\in\DD$, let
\[
D_{\Gamma,\oo}=\{z\in\DD; d(z,\oo) <
d(z,\gamma(\oo)),\quad\forall\,\gamma\in\Gamma\}
\]
\noindent Denote the corresponding Dirichlet domain, a convex
fundamental domain with compact closure in $\DD$, admitting an
even number of geodesic sides and an even number of vertices, some
of which may be elliptic. More precisely, the boundary of
$D_{\Gamma,\oo}$ is  a disjoint union of semi-closed geodesic
segments $S^L_{-r},\cdots,S^L_{-1},S^L_1,\cdots,S^L_{r}$, closed
to the left and open to the right, or, equivalently, to a union of
semi-closed geodesic segments
$S^R_{-r},\cdots,S^R_{-1},S^R_1,\cdots,S^R_{r}$, closed to the
right and open to the left;  for each $k$, the intervals  $S^L_k$
and $S^R_k$ have the same endpoints and $S^L_k$ is associated to
$S^R_{-k}$ by an element $a_k\in\Gamma$ satisfying
$a_k(S^L_k)=S^R_{-k}$. The elements $a_k$ generate $\Gamma$ and
satisfy $a_{-k}=a_k^{-1}$, for $k=\pm 1,\cdots,\pm r$.

To define the two Bowen-Series transformations $T_L$ and $T_R$
geometrically, we need to impose a geometric condition on
$\Gamma$: following \cite{BoSe1}, \cite{Ser2} and \cite{Ser1}, we
say that $\Gamma$ satisfies the {\it even corner} property if, for
each $1\leq|k|\leq r$, the complete geodesic line through $S^L_k$
is equal to a disjoint union of $\Gamma$-translates of the sides
$S^L_l$, with $ 1\leq|l|\leq r$. Some $\Gamma$ do not satisfy this
geometric property. Nevertheless, any two co-compact Fuchsian
groups $\Gamma$ and $\Gamma'$, with identical signature, are
geometrically isomorphic, that is, there exists a group
isomorphism $h_*:\Gamma\to\Gamma'$ and a quasi-conformal
orientation preserving homeomorphism $h:\DD\to\DD$ admitting an
extension to a conjugating homeomorphism
$h:\partial\DD\to\partial\DD$, that is,
\[
h(\gamma(z))=h_*(\gamma)(h(z)),\quad\forall\ \gamma\in\Gamma.
\]

An important observation in \cite{BoSe1}, \cite{Ser2} and
\cite{Ser1} is that any co-compact Fuchsian group is geometrically
isomorphic to a Fuchsian group with identical signature and
satisfying the {\it even corner} property. We are going to recall
the Bowen and Series construction in the case  that $\Gamma$
possesses the {\it even corner} property and will show that their
main conclusions remain valid under geometric isomorphisms.

The complete geodesic line associated to a side $S_k^L$ cuts the
boundary at infinity $\SS^1$ at two points $s_k^L$ and $s_k^R$,
positively oriented with respect to $s_k^L$, the oriented geodesic
line $]]s_k^L,s_k^R[[$ seeing the origin $\oo$ to the left. Both
end points $s_k^L$ and $s_k^R$ are  neutrally stable with respect
to the associated generator $a_k$, that is,
$|a_k'(s_k^L)|=|a_k'(s_k^R)|=1$. The family of open intervals
$]s_k^L,s_k^R[$ covers $\SS^1$; since these intervals
$]s_k^L,s_k^R[$ overlap each other, there is no canonical
partition adapted to this covering. Nevertheless, we may associate
two well defined partitions, the left partition $\aaa_L$ and  the
right partition $\aaa_R$. The former consists of disjoint
half-closed intervals,
\[
\aaa_L=\{A_{-r}^L,\cdots,A_{-1}^L,A_1^L,\cdots,A_{r}^L\},
\]
\noindent given by $A_k^L=[s_k^L,s_{l(k)}^L[$ where $s_{l(k)}^L$
denotes the nearest point $s_l^L$ after $s_k^L$, according to a
positive orientation. Each $A_k^L$ belongs to the unstable domain
of the hyperbolic element $a_k$, that is, $|a_k'(\xi)|\geq 1$, for
each $\xi\in A_k^L$. By definition, the left Bowen-Series
transformation $T_L:\SS^1\mapsto\SS^1$ is given by
\[
T_L(\xi)=a_k(\xi),\qquad\text{if}\quad \xi\in A_k^L.
\]
\noindent Analogously, $\SS^1$ can be partitioned into half-closed
intervals
\[
\aaa_R=\{A_{-r}^R,\cdots,A_{-1}^R, A_1^R,\cdots,A_{r}^R\},
\]
\noindent where $A_k^R=]s_{j(k)}^R,s_k^R]$, and $s_{j(k)}^R$
denotes the nearest $s_j^R$ before $s_k^R$, according to a
positive orientation. The right Bowen-Series transformation is
given by
\[
T_R(\eta)=a_k(\eta),\qquad\text{if}\quad \eta\in A_k^R.
\]
\noindent The two partitions $A^L $ and $A^R$ generate two ways of
coding a trajectory. Let $\gamma_L:\SS^1\mapsto\Gamma$ and
$\gamma_R:\SS^1\mapsto\Gamma$ be the left and right symbolic
coding defined by
\[
\gamma_L[\xi]=a_k,\quad\text{if $\xi\in A_k^L$}, \qquad
\text{and}\,\,\,\,\gamma_R[\eta]=a_k,\quad\text{if $\eta\in
A_k^R$}.
\]
\noindent In particular, $T_R(\eta)=\gamma_R[\eta](\eta)$ and
$T_L(\xi)=\gamma_L[\xi](\xi)$, for each $\xi\in\SS^1$. Also, it is
known  that $T_R^2$ and $T_L^2$ are expanding. Series, in
\cite{Ser2}, \cite{Ser3} and \cite{Ser1}, and later, Adler and
Flatto in \cite{AdFl}, proved that $T_L$ (respectively $T_R$) is
Markov with respect to a partition of $\ii^L=\{I^L_k\}_{k=1}^q$
(respectively $\ii^R=\{I^R_l\}_{l=1}^q$) that is  finer than
$\aaa_L$ (respectively $\aaa_R$). The semi-closed intervals
$I^L_k$ and $I^R_l$  are of the same kind as $A_k^L$ and $A_l^R$,
and have the same closure.

\begin{definition}
A dynamical system $(\SS^1,T,\{I_k\})$ is said to be a piecewise
$\Gamma$-M\"obius Markov transformation if $T:\SS^1\to\SS^1$ is a
surjective map, and  $\{I_k\}$ is a finite partition of $\SS^1$
into intervals such that:
\begin{enumerate}
\item for each $k$, $T(I_k)$ is a union of adjacent intervals
$I_l$; \item for each $k$, the restriction of $T$ to $I_k$
coincides with an element $\gamma_k\in\Gamma$; \item some finite
iterate of $T$ is uniformly expanding.
\end{enumerate}
\end{definition}

\begin{theorem}{\rm(\cite{BoSe1}, \cite{Ser1})} For any co-compact Fuchsian group $\Gamma$, there exists a piecewise $\Gamma$-M\"obius Markov transformation $(\SS^1,T,\{I_k\})$ which is transitive and orbit equivalent to $\Gamma$.
\end{theorem}

The Ruelle transfer operator can be defined for any picewise
$\cc^2$ Markov transformation  $(\SS^1,T,\{I_k\})$ and any
potential function $A$. Actually, we need a particular complex
transfer operator given by the potential
\[
A=-s\ln|T'|.
\]
\noindent For any function $\psi:\SS^1\to\CC$, define
\begin{gather*}
(\ll_s(\psi))(\xi') = \sum_{T(\xi)=\xi'}e^{A(\xi)}\psi(\xi) =
\sum_{T(\xi)=\xi'}\frac{\psi(\xi)}{|T'(\xi)|^s},
\end{gather*}
\noindent where the summation is taken other all preimages $\xi$
of $\xi'$ under $T$. We modify $\ll_s$ slightly,  so that it acts
on the space of piecewise $\cc^1$ functions. Let $\{I_k\}_{k=1}^q$
be a partition of $S^1$. Given a piecewise $C^1$ function and
$\oplus_{k=1}^q\psi_k\in\oplus_{k=1}^q\cc^1(\bar I_k)$ set
\[
\ll_s^L\psi=\oplus_{l=1}^q\phi_l, \quad \textrm{where}\quad
\phi_l=\sum_{I_l\subset T(I_k)} \frac{\psi_k\circ
T_{k,l}^{-1}}{|T'\circ T_{k,l}^{-1}|^s},
\]
\noindent and  $T_{k,l}^{-1}$ denotes  the restriction to $I_l$
of the inverse of $T:I_k\to T(I_k) \supset I_l$.

\begin{proposition}\label{label4} Let $\Gamma$ be a co-compact Fuchsian group. Let $s=\frac{1}{2}+it$ and
$f$ be an automorphic eigenfunction of $-\Delta$, that is, $\Delta
f=-s(1-s)f$. Let $(\SS^1,T,\{I_k\})$ be a piecewise
$\Gamma$-M\"obius Markov transformation and $\ll_s$ be the Ruelle
transfer operator corresponding to the observable $A=-s\ln|T'|$.
Then the Helgason distribution $\dd_{f,s}$ satisfies
\[
(\ll_s)^*\dd_{f,s}=\dd_{f,s}.
\]
\end{proposition}

\begin{proof}
Let $\oplus_{k=1}^q\psi_k$ be a piecewise $\cc^1$ function in
$\oplus_{k=1}^q\cc^1(\bar I_k)$. Using Proposition \ref{label1},

\begin{align*}
\int(\ll_s\psi)(\xi)\,\dd_{f,s}(\xi)
&= \sum_{l=1}^q \int_{I_l}(\ll_s\psi)_l(\xi)\,\dd_{f,s}(\xi) \\
&=  \sum_{T(I_k)\supset I_l} \int_{I_l} \frac{\psi_k\circ T_{k,l}^{-1}}{|T'\circ T_{k,l}^{-1}|^s}(\xi)\,\dd_{f,s}(\xi) \\
&= \sum_{T(I_k)\supset I_l} \int_{T^{-1}(I_l)\cap I_k} \psi_k(\xi)\,\dd_{f,s}(\xi) \\
&= \sum_{k=1}^q \int_{I_k} \psi_k(\xi)\,\dd_{f,s}(\xi) = \int \psi(\xi)\,\dd_{f,s}(\xi).
\end{align*}
\end{proof}

Series in \cite{Ser1}, Adler and Flatto in \cite{AdFl},  and
Morita in \cite{Mo}  noticed that $T_L$ admits a natural extension
$\hat T:\hat \Sigma\mapsto\hat\Sigma$ strongly related to $T_R$.
We also showed the existence of such a $\hat T$ in \cite{LT}, and
it was an important step in the proof of Theorem 3 of \cite{LT}.
The following definition explains how the two maps $T_L$ and $T_R$
are glued together in an abstract way.

\begin{definition}
Let $\Gamma$ be a co-compact Fuchsian group. A dynamical system
$(\hat\Sigma,\hat T,\{I_k^L\},\{I_l^R\},J)$ is said to be a
piecewise $\Gamma$-M\"obius baker transformation if it admits a
description as follows.
\begin{enumerate}
\item $\{I_k^L\}$ and $\{I_l^R\}$ are finite partitions of $\SS^1$
into disjoint intervals; $J(k,l)$ is a $\{0,1\}$-valued function,
and $\hat\Sigma$ is the subset of $\SS^1\times\SS^1$ defined by
\[
\hat\Sigma=\coprod_{J(k,l)=1}I_k^L\times I_l^R.
\]
\item For each $k$, $Q_k^R=\coprod\{I_l^R; J(k,l)=1\}$ is an
interval whose closure is disjoint from $\bar I_k^L$. For each
$l$, $Q_l^L=\coprod\{I_k^L; J(k,l)=1\}$ is an interval whose
closure is disjoint from $\bar I_l^R$. Let $I^L(\xi)=I^L_k$ and
$Q^R(\xi)=Q^R_k$, for $\xi\in I^L_k$. Let $I^R(\eta)=I^R_l$ and
$Q^L(\eta)=Q^L_l$, for $\eta\in I^R_l$. \item $\hat
T:\hat\Sigma\to\hat \Sigma$ is bijective and is given by
\[
\left\{\begin{array}{ccc}
\hat T(\xi,\eta) &=& (T_L(\xi),S_R(\xi,\eta)), \\
\hat T^{-1}(\xi',\eta') &=& (S_L(\xi',\eta'),T_R(\eta')),
\end{array}\right.
\]
\noindent for certain maps $T_L,T_R:\SS^1\to \SS^1$ and
$S_L,S_R:\hat\Sigma\to \SS^1$. \item $(\SS^1,T_L,\{I_k^L\})$ and
$(\SS^1,T_R,\{I_l^R\})$ are piecewise $\Gamma$-M\"obius Markov
transformations. There exist two functions
$\gamma_L:\SS^1\to\Gamma$, respectively $\gamma_R:\SS^1\to\Gamma$,
that are piecewise constant on each $I_k^L$, respectively
$\{I_l^R\}$, and satisfying
\[
\left\{\begin{array}{ccc}
\hat T(\xi,\eta) &=& (\gamma_L[\xi](\xi),\gamma_L[\xi](\eta)) \\
\hat T^{-1}(\xi',\eta') &=& (\gamma_R[\eta'](\xi'),\gamma_R[\eta'](\eta'))
\end{array}\right.
\]
\end{enumerate}

\noindent The maps $T_L$ and $T_R$ are called the left and right
Bowen-Series transformations., whereas $\gamma_L$ and $\gamma_R$
are the left and right Bowen-Series codings. Finally, we say that
$J$ is the incidence matrix, which we extend as a function on
$\SS^1\times\SS^1$ defining
\[
\left\{\begin{array}{l}
J(\xi,\eta)=1, \quad\textrm{if}\quad (\xi,\eta)\in\hat\Sigma ,\\
J(\xi,\eta)=0, \quad\textrm{if}\quad (\xi,\eta)\not\in\hat\Sigma.
\end{array}\right.
\]
\end{definition}

Notice that this definition is equivariant by geometric
isomorphisms. For co-compact Fuchsian groups satisfying the {\it
even corner} property, Adler and Flatto in \cite{AdFl}, Series in
\cite{Ser1} (and, for a particular example, in \cite{LT}) obtained
geometrically the existence of a piecewise $\Gamma$-M\"obius baker
transformation with left $T_L$ and right $T_R$  maps orbit
equivalent to $\Gamma$. By geometric isomorphism considerations,
we obtain more generally the following.

\begin{proposition}{\rm (\cite{AdFl}, \cite{Ser1}, \cite{LT})}\label{label7}
For any co-compact Fuchsian group $\Gamma$, there exists a
picewise $\Gamma$-M\"obius baker transformation with left and
right Bowen-Series transformations that are transitive and orbit
equivalent to $\Gamma$.
\end{proposition}

The two maps $T_L$ and $T_R$ are related to the action of the
group $\Gamma$ on the boundary $\SS^1$. The baker transformation
$(\hat \Sigma,\hat T)$ encodes this action into a unique dynamical
system. For later reference, we state two further properties of
this encoding.
\begin{remark}\quad\par\rm
\begin{enumerate}
\item  The two codings $\gamma_L$ and $\gamma_R$ are reciprocal,
in the following sense:
\[
\gamma_R[\eta']=\gamma_L^{-1}[\xi],
\quad\textrm{whenever}\quad (\xi',\eta')=\hat T(\xi,\eta).
\]
\item For any $\xi'$ and $\eta$ in $\SS^1$, there is a bijection between the two finite sets
\[
\{\xi;(\xi,\eta)\in\hat\Sigma \textrm{ and } T_L(\xi) = \xi'\},
\quad
\{\eta';(\xi',\eta')\in\hat\Sigma \textrm{ and }  T_R(\eta')=\eta\}.
\]
\end{enumerate}
\end{remark}

In order to better understand this baker transformation, we briefly explain how $(\hat\Sigma,\hat T)$ is conjugate to a specific Poincar\'e section of the geodesic flow on the surface $N=T^1M$. We assume for the rest of this section that $\Gamma$ satisfies the {\it even corner} property.

 Since $D_{\Gamma,\oo}$ is a convex fundamental domain, every geodesic (modulo $\Gamma$) cuts $\partial D_{\Gamma,\oo}$ at two distinct points $p$ and $q$, unless the geodesic is tangent to one of the sides of $D_{\Gamma,\oo}$. These tangent geodesics correspond to a finite union of closed geodesics. We could have parametrized the set of oriented geodesics by all pairs $(p,q)\in\partial D_{\Gamma,\oo}\times\partial D_{\Gamma,\oo}$, with $p$ and $q$ not belonging to the same side of $D_{\Gamma,\oo}$, but we prefer to introduce the space $X$ of all $(x,y)\in\SS^1\times \SS^1$ oriented geodesics  $[[y,x]]$, either cutting the interior of $D_{\Gamma,\oo}$ or passing through one of the corners of $D_{\Gamma,\oo}$ and seeing $\oo$ to the left. Using these notations, we define  the two intersection points $p=p(x,y)\in\partial D_{\Gamma,\oo}$ and $q=q(x,y)\in\partial D_{\Gamma,\oo}$ for every oriented geodesic $[[y,x]]$, $(x,y)\in X$, such that $[[q,p]]=[[y,x]]\cap\bar D_{\Gamma,\oo}$ has the same orientation as $[[y,x]]$.

For a geodesic passing through a corner, $p=q$, unless the geodesic is tangent to a side of $D_{\Gamma,\oo}$. We are now in a position to define a geometric Poincar\'e section $B:X\to X$. If $(x,y)\in X$, the geodesic $[[y,x]]$ leaves $D_{\Gamma,\oo}$ at $p=p(x,y)\in S_i$, for some side $S^L_i$. Since $S^L_i$ and $S^R_{-i}$ are permuted by the generator $a_i$, the new geodesic $a_i([[y,x]])=[[y',x']]$ enters again the fundamental domain at a new point $q'=q(x',y')$ with $q'=a_i(p)\in S^R_{-i}$. By definition, $B(x,y)=(x',y')$ and the map $B:X\to X$ is called a geodesic billiard like the codings as for $T_L$ and $T_R$, we introduce two geometric codings $\gamma_B:X\to\Gamma$ and $\bar\gamma_B:X\to\Gamma$ given by
\[
\left\{\begin{array}{ccc}
\gamma_B[x,y]=a_i &\textrm{if}& p(x,y)\in S^L_i, \\
\bar \gamma_B[x,y]=a_i &\textrm{if}& q(x,y)\in S^R_i.
\end{array}\right.
\]
\noindent Now the geodesic billiard can  be defined by
\[
\left\{\begin{array}{ccc}
B(x,y) &=& (\gamma_B[x,y](x),\gamma_B[x,y](y)), \\
B^{-1}(x',y') &=& (\bar\gamma_B[x',y'](x'),\bar\gamma_B[x',y'](y'))
.\end{array}\right.
\]
\noindent Notice that $\bar \gamma_B\circ B=\gamma_B^{-1}$. The map $B$ is very close to be a  baker transformation: $B$ and $B^{-1}$ have the same structure as $\hat T$ and $\hat T^{-1}$, and $\gamma_B$ (respectively, $\bar \gamma_B$) plays the role as $\gamma_L$ (respectively, $\gamma_R$). The main difference is that $\gamma_B[x,y]$ depends on both $x$ and $y$, but $\gamma_L[\xi]$ depends only on $\xi$. Nevertheless, we have the following crucial result.

\begin{theorem}{\rm(\cite{AdFl}, \cite{Ser1}, \cite{LT})}\label{label9} There exists a $\Gamma$-M\"obius baker transformation $(\hat\Sigma,\hat T)$ conjugate to $(X,B)$. More precisely, there exists a map $\rho:X\to \Gamma$ such that $\pi(x,y)=(\rho[x,y](x),\rho[x,y](y))$,  defines a conjugating map $\pi:X\to \hat \Sigma$ between $\hat T$ and $B$, such that $\hat T\circ \pi=\pi\circ B$. Equivalently, $\gamma_L\circ\pi$ and $\gamma_B$ are cohomologous over $(X,B)$, that is, $\gamma_L\circ\pi\rho=\rho\circ B\gamma_B$, and $\gamma_R\circ \pi$ and $\bar \gamma_B$ are cohomologous over $(X,B)$, that is, $\gamma_R\circ\pi \rho=\rho\circ B^{-1}\bar\gamma_B$.
\end{theorem}

\section{Proof of Theorem \ref{label0}}

We want to associate to any eigenfunction $f$ of the Laplace operator a nonzero piecewise real analytic function $\psi_{f,s}$ that is a solution of the functional equation
\[
\ll_s^L(\psi_{f,s})=\psi_{f,s},\quad
\text{where}\quad
\ll_s^L(\psi)(\xi')=\sum_{T_L(\xi)=\xi'}\frac{\psi(\xi)}{|T'_L(\xi)|^s}.
\]
\noindent The main idea is to use a kernel $k(\xi,\eta)$ introduced in Theorem 7 of \cite{BLT}, as well by in Haydn in \cite{Hay}, and by Bogomolny and Carioli in \cite{BC} and  \cite{BC2}, in the context of double-sided subshifts of finite type. We begin by extending this definition to include  baker transformations.

\begin{definition}
Let $(\hat\Sigma,\hat T)$ be a piecewise $\Gamma$-M\"obius baker transformation, with $T_L$ and $T_R$ the left and right Bowen-Series transformations. Let $A_L:\SS^1\to\CC$ and $A_R:\SS^1\to\CC$ be two potential functions. We say that $A_L$ and $A_R$ are in involution if there exists a nonzero kernel $k:\hat\Sigma\to\CC^*$, called an involution kernel, such that
\[
k(\xi,\eta)e^{A_L(\xi)}=k(\xi',\eta')e^{A_R(\eta')},\quad
\textrm{whenever}
\quad
(\xi',\eta')=\hat T(\xi,\eta)  \in\hat\Sigma.
\]
\noindent The kernel $k$ is extended to $\SS^1\times\SS^1$ by $k(\xi,\eta)=0$, for $(\xi,\eta)\not\in\hat\Sigma$.
\end{definition}

\begin{remark}\quad\par\rm
\begin{enumerate}
\item Let $W(\xi,\eta)=\ln k(\xi,\eta)$, for  $(\xi,\eta)\in\hat\Sigma$. Then $A_L$ and $A_R$ are cohomologous, that is $A_L-A_R\circ\hat T=W\circ \hat T-W$.
\item If $A_L(\xi)$ is H\"older, then there exists a H\"older function $A_R(\eta)$ (depending only on $\eta$) in involution with $A_L$ with a H\"older involution kernel.
\item If $\ll_L$ and $\ll_R$ are the two Ruelle transfer operators associated to $A_L$ and $A_R$, if $A_L$ and $A_R$ are in involution with respect to a kernel $k$, and if $\nu$ is an eigenmeasure of $\ll_R$, that is, $\ll_R^*(\nu)=\lambda\nu$, then $\psi(\xi)=\int k(\xi,\eta)\,d\nu(\eta)$ is an eigenfunction of $\ll_L$,
    that is, $\ll_L(\psi)=\lambda\psi$.
\end{enumerate}
\end{remark}

These remarks appeared  first in \cite{Hay} and were later rediscovered  in \cite{BLT}, in the context of a subshift of finite type. The proofs in this general context can be easily reproduced. The third remark suggests a strategy to obtain the eigenfunction $\psi_{f,s}$, by taking $A_L=-s\ln|T'_L|$, $A_R=-s\ln|T'_R|$ and replacing $\nu$ by the distribution $\dd_{f,s}$. All there is left to prove is that $-\ln|T'_L|$ and $-\ln|T'_R|$ are in involution with respect to a piecewise $\cc^1$ involution kernel. It so happens that this  involution kernel exists and is given by the Gromov distance.

\begin{definition}
The Gromov distance  $d(\xi,\eta)$ between two points  $\xi $ and $\eta $ at infinity is given by
\[
d^2(\xi,\eta)=\exp\Big(
-b_\xi(\oo,z)-b_\eta(\oo,z)
\Big),
\]
\noindent for any point $z$ on the geodesic line $[[\xi,\eta]]$. Notice that this definition depends on the choice of the origin $\oo$ (but not on $z\in[[\xi,\eta]]$).
\end{definition}

In the Poincar\'e disk model, $(\xi,\eta)\in\SS^1\times\SS^1$, or in the upper half-plane, $(s,t)\in\RR\times\RR$, the Gromov distance takes the simple form
\[
d^2(\xi,\eta)=\frac{1}{4}|\xi-\eta|^2,
\quad\textrm{or}\quad
d^2(s,t)=\frac{|s-t|^2}{(1+s^2)(1+t^2)}.
\]

\begin{lemma}\label{label8}
Let $T_L:\SS^1\to\SS^1$ and $T_R:\SS^1\to\SS^1$ be the two left and right Bowen-Series transformations of a $\Gamma$-M\"obius Markov baker transformation  $(\hat\Sigma,\hat T)$. Then the two potential functions $A_L(\xi)=-\ln|T'_L(\xi)|$ and $A_R(\eta)=-\ln|T'_R(\eta)|$ are in involution and
\[
A_L(\xi)-A_R(\eta')=W(\xi',\eta')-W(\xi,\eta),\quad
\textrm{for}
\quad
(\xi',\eta')=\hat T(\xi,\eta)\in\hat\Sigma,
\]
\noindent where $W(\xi,\eta)=b_\xi(\oo,z)+b_\eta(\oo,z)$ and $z$ is any point of the geodesic line $[[\xi,\eta]]$. In particular, $k(\xi,\eta)=\exp(W(\xi,\eta))=4/d^2(\xi,\eta)$ is an involution kernel.
\end{lemma}

\begin{proof}[ of Lemma \ref{label8}.] To simplify the notation, we call $(\xi',\eta')=\hat T(\xi,\eta)$, $\gamma_L=\gamma_L[\xi]$, and $\gamma_R=\gamma_R[\eta']$. We also recall the relation $\gamma_R=\gamma_L^{-1}$. Then, choosing any point $z\in[[\xi,\eta]]$, we get
\begin{align*}
A_L(\xi)-A_R(\eta') &=
-b_\xi(\oo,\gamma_L^{-1}\oo)+b_{\eta'}(\oo,\gamma_R^{-1}\oo) \\
&= -b_\xi(\oo,z)-b_\xi(z,\gamma_L^{-1}\oo) \\
&\quad+ b_{\eta'}(\oo,\gamma_L(z))+b_{\eta'}(\gamma_L(z),\gamma_R^{-1}\oo) \\
&= W(\xi',\eta')-W(\xi,\eta),
\end{align*}
\noindent where $W(\xi',\eta')=b_{\eta'}(\oo,\gamma_L(z))-b_\xi(z,\gamma_L^{-1}\oo)$ and $W(\xi,\eta)=b_\xi(\oo,z)-b_{\eta'}(\gamma_L(z),\gamma_R^{-1}\oo)$.
\end{proof}

\noindent Notice that if $A(\xi)$ and $\bar A(\eta)$ are in involution by a positive kernel $k(\xi,\eta)$, then $sA(\xi)$ and $s\bar A(\eta)$ are in involution by $k(\xi,\eta)^s$.

\begin{lemma}\label{label10}
Let $T_L:\SS^1\to\SS^1$ and $T_R:\SS^1\to\SS^1$ be the two left and right Bowen-Series transformations of a $\Gamma$-M\"obius Markov baker transformation  $(\hat\Sigma,\hat T)$. Let $A_L:\SS^1\to \RR$ and $A_R:\SS^1\to\RR$ be two potential functions in involution with respect to a kernel $k(\xi,\eta)$. Let $\ll_L$ and $\ll_R$ be the two Ruelle transfer operators associated to $A_L$ and $A_R$. Then, for any $\xi'\in\SS^1$ and $\eta\in\SS^1$,
\[
\ll_R(k(\xi',\cdot))(\eta)=\ll_L(k(\cdot,\eta))(\xi').
\]
\end{lemma}

\begin{proof} Given $\xi'\in\SS^1$ and $\eta\in\SS^1$,  the two finite sets
\[
\{\eta'\in\SS^1;\ T_R(\eta')=\eta,\,\, J(\xi',\eta')=1\},
\quad
\{\xi\in\SS^1;\ T_L(\xi)=\xi',\,\, J(\xi,\eta)=1\}
\]
\noindent are in bijection. Thus, we obtain
\begin{align*}
\ll_R(k(\xi',\cdot))(\eta) &= \sum_{T_R(\eta') = \eta} k(\xi',\eta')e^{A_R(\eta')} \\
&= \sum_{T_L(\xi)=\xi'} k(\xi,\eta)e^{A_L(\xi)}=
\ll_L(k(\cdot,\eta))(\xi')
\end{align*}
\end{proof}

Theorem \ref{label0} now follows immediately from lemmas \ref{label8} and \ref{label10}.

\begin{proof}[ of Theorem \ref{label0}.] We first prove that $\psi_{f,s}(\xi)=\int k(\xi,\eta)^s\dd_{f,s}(\eta), $ with $k(\xi,\eta)=J(\xi,\eta)/d^2(\xi,\eta)$, is a solution of the equation $\ll_s^L\psi_f=\psi_f$. In fact, we have
\begin{align*}
\psi_{f,s}(\xi') &= \int k^s(\xi',\eta')\,\dd_{f,s}(\eta')
= \int \ll_s^R(k^s(\xi',\cdot))(\eta)\,\dd_{f,s}(\eta) \\
&= \int \ll_s^L(k^s(\cdot,\eta)(\xi')\,\dd_{f,s}(\eta)
= (\ll_s^L \psi_{f,s})(\xi').
\end{align*}
\noindent We next prove that $\psi_{f,s}\not=0$. Suppose on the contrary that $\psi_{f,s}(\xi')=0$ for each $\xi'\in\SS^1$. Following Haydn \cite{Hay}, we introduce step functions of the form
\[
\bar \chi(\xi',\eta')=\chi\circ pr_1\circ \hat T^{-1}(\xi',\eta'),
\]
\noindent where $\chi=\chi(\xi)$ depends only on $\xi$. For instance, for some fixed $\xi'$, let $\chi$ be the characteristic function of the interval $I^L(n,\xi)=\cap_{k=0}^n T_L^{-k}(I^L\circ T_L^k(\xi)),$ for some $\xi$ such that $T_L^n(\xi)=\xi'$. Let $Q^R(\xi)=\{\eta\in\SS^1 ; J(\xi,\eta)=1\}$ and write
\[
\gamma_L[n,\xi]=\gamma_L[T^{n-1}_L(\xi)]\cdots\gamma_L[T_L(\xi)]\gamma_L[\xi],\quad
Q^R(n,\xi)=\gamma_L[n,\xi]Q^R(\xi).
\]
\noindent Then $\bar \chi$ equals the characteristic function of the rectangle $I^L(\xi')\times Q^R(n,\xi)$ and $Q^R(\xi')$ is equal to the disjoint union of the intervals $Q^R(n,\xi)$, for all $\xi$ such that $T^n_L(\xi)=\xi'$. We also denote by $\Delta(\xi')$ the set of endpoints of $Q^R(n,\xi)$, for all $T_L^n(\xi)=\xi'$, and observe that $\Delta(\xi')$ is a dense subset of $Q^R(\xi')$. Using the same ideas as in Lemma \ref{label10}, we obtain
\[
\int \bar\chi(\xi',\eta')k^s(\xi',\eta')\dd_{f,s}(\eta') =
(\ll_s^L)^n(\chi\psi_{f,s})(\xi')=0,\quad
\forall\ \xi'\in\SS^1.
\]
\noindent In particular, if $\tilde\alpha(\xi')<\tilde\beta(\xi')<\tilde\alpha(\xi')+2\pi$ are chosen such that $\exp i\tilde\alpha(\xi')$ and $\exp i\tilde\beta(\xi')$ are the two endpoints of the interval $Q^R(\xi')$, if $\tilde k(\theta)=k(\xi',\exp i\theta)$, then
\[
\tilde k(\beta)\tilde\dd_{f,s}(\beta) = \tilde k(\tilde\alpha(\xi'))\tilde\dd_{f,s}(\tilde\alpha(\xi')) + \int_{\tilde\alpha(\xi')}^\beta \frac{\partial\tilde k}{\partial\theta}\tilde\dd_{f,s}(\theta)\,d\theta.
\]
\noindent for every $\beta\in[\tilde\alpha(\xi'),\tilde\beta(\xi')]\cap\Delta(\xi')$. Since $\tilde k(\theta)\not=0$, for  each $\theta\in[\tilde\alpha(\xi'),\tilde\beta(\xi')]$, we conclude  that the above equality applies to all $\beta\in[\tilde\alpha(\xi'),\tilde\beta(\xi')]$, the two functions $\tilde k(\beta)\tilde\dd_{f,s}(\beta)$ and $\tilde\dd_{f,s}(\beta)$ are $\cc^1$, and
\[
\int_{\tilde\alpha(\xi'))}^\beta k(\theta) \frac{\partial \tilde\dd_{f,s}}{\partial\theta}\,d\theta=0,\quad
\forall\ \beta\in [\tilde\alpha(\xi'),\tilde\beta(\xi')].
\]
\noindent Therefore, $\tilde\dd_{f,s}(\theta)$ is a constant function on each $[\tilde\alpha(\xi'),\tilde\beta(\xi')]$, thus everywhere on $\SS^1$. It follows that the distribution $\dd_{f,s}$ would have to be equal to zero, which is impossible, because it represents a nonzero eigenfunction $f$.

\end{proof}

We would like to thank I. Efrat for showing us the reference
\cite{BC}, F. Ledrappier for the references \cite{Le}  and
\cite{LeZa} and C. Doering for help us with corrections in our
text. Finally, we would like to thank the referees for their
careful reading and comments.

\vspace{0.9cm}

E-mail:
\vspace{0.2cm}
arturoscar.lopes@gmail.com
\par
philippe.thieullen@math.u-bordeaux1.fr

\bigskip
AMS-MSC : 37C30, 11F12, 11F72, 46F12


\begin{thebibliography}{999}


\bibitem{AdFl} R. Adler, L, Flatto, Geodesic flows, interval maps and symbolic dynamics, {\it Bull. Amer. Math. Soc.}, vol. 25 (1991), pp. 229--334.

\bibitem{AZ} N. Anantharaman and S. Zelditch, Patterson-Sullivan distributions and quantum ergodicity, {\it Ann. Henri Poincar\'e}, vol. 8 (2007), pp. 361--426.

\bibitem{BLT} A. Baraviera,  A. O. Lopes and Ph. Thieullen,
A large deviation principle for  equilibrium states of H\"older
potentials: the zero temperature case, {\it Stochastics and
Dynamics}, vol. 6:1 (2006).

\bibitem{BeKeSe} T. Bedford, M. Keane and C. Series, {\it Ergodic Theory, Symbolic Dynamics and Hyperbolic Spaces}, Oxford Univ Press (1991).

\bibitem{BeMa} M. Bekka and M. Mayer, {\it Ergodic Theory and Topological Dynamics of Group Actions on Homogeneous Spaces}, Cambridge Univ. Press (2000).


\bibitem{BC} E.B. Bogomolny and M. Carioli, Quantum maps of geodesic flows on surfaces of constant negative curvature, {\it IV International Conference on ``Path Integrals from meV to MeV"}, Tutzing May,  (1992), pp. 18--21.

\bibitem{BC2} E.B. Bogomolny and M. Carioli, Quantum maps from transfer operator. {\it Physica D}, vol. 67 (1993), pp. 88--112.

\bibitem{BoSe1} R. Bowen and C. Series, Markov maps associated to Fuchsian groups, {\it Pub. Math. IHES}, vol. 50 (1979), pp. 153-170.

\bibitem{FF} L. Flaminio and G. Forni, Invariant distributions and the time averages for horocycle flows, {\it Duke Math. J.}, vol. 119 (2003), pp. 465--526.

\bibitem{Hay} N. Haydn. Gibbs' functionals on subshifts. {\it Commun. Math. Phys.} vol. 134 (1990), pp. 217--236.

\bibitem{He1} S. Helgason, {\it Analysis on Lie Groups and Homogeneous Spaces}, A.M.S., 1972.

\bibitem{He} S. Helgason, {\it Topics in Harmonic Analysis on Homogeneous
Spaces}, Birkhauser, 1981.

\bibitem{I} H. Iwaniec, {\it Spectral Methods of Automorphic Forms}, AMS (2002).

\bibitem{Le} J. B. Lewis, Spaces of holomorphic functions equivalent to the even Maass cusp forms, {\it Invent. Math.}, vol. 127 (1997), pp. 271--306.

\bibitem{LeZa} J. Lewis, D. Zagier, Period functions for Maass wave forms, {\it Annals of Mathematics}, vol. 153 (2001), pp. 191--258.

\bibitem{LT} A.O. Lopes  and P. Thieullen, Mather theory and
the Bowen-Series transformation, {\it Annal. Inst. Henry Poincar\'e,
Analyse non-lin\'eaire},  vol. 23:5 (2006), pp. 663--682.

\bibitem{Ma} D. H. Mayer,
Thermodynamic formalism approach to Selberg's zeta function for
PSL $(2,\mathbb{Z})$, {\it Bull. Amer. Math. Soc.}, vol. 25:1 (1991),
pp. 55--60.

\bibitem{Mo} T. Morita, Markov systems and transfer operators associated with cofinite Fuchsian groups. {\it Erg. Theory and Dyn. Syst.} vol. 17 (1997), pp. 1147--1181.

\bibitem{Ot} J-P. Otal,
Sur les fonctions propres du laplacien du disque hyperbolique, {\it C.
R. Acad. Sci. Paris S\'er. I Math.}, vol. 327:2 (1998), pp. 161--166.

\bibitem{Pa} S.J. Patterson, The limit set of a Fuchsian group, {\it Acta Math.}, vol. 136 (1976), pp. 241--273.

\bibitem{Po} M. Pollicott, Some applications of thermodynamic formalism to manifolds with constant negative curvature, {\it Advances in Mathematics}, vol. 85 (1991), pp. 161--192.

\bibitem{Ser2} C. Series, Symbolic dynamics for geodesic flows. {\it Acta Math.}, vol. 146 (1981), pp. 103--128.

\bibitem{Ser3} C. Series, The modular surface and continuous fractions, {\it J. London Math. Soc.}, vol. 31 (1985), pp. 69--80.

\bibitem{Ser1} C. Series, Geometrical Markov coding on surfaces of constant negative curvature,
{\it Erg. Th. and Dynam. Sys.}, vol. 6 (1986), pp. 601--625.

\bibitem{Su} D. Sullivan, The density at infinity of a discrete group of hyperbolic motions, {\it Publ. Math. IHES}, vol. 50 (1979), pp. 171--202.

\bibitem{Ze} S. Zelditch, Uniform distribution of eigenfunctions on compact hyperbolic surfaces, {\it Duke Math. J.}, vol. 55 (1987), pp. 919--941.

\end{thebibliography}
\end{document}